# Estimation of population-level summaries in general semiparametric repeated measures regression models[*]

**Arnab Maity**[1], **Tatiyana V. Apanasovich**[2] and **Raymond J. Carroll**[1]

*Texas A&M University, Cornell University and Texas A&M University*

**Abstract:** This paper considers a wide family of semiparametric repeated measures regression models, in which the main interest is on estimating population-level quantities such as mean, variance, probabilities etc. Examples of our framework include generalized linear models for clustered/longitudinal data, among many others. We derive plug-in kernel-based estimators of the population level quantities and derive their asymptotic distribution. An example involving estimation of the survival function of hemoglobin measures in the Kenya hemoglobin study data is presented to demonstrate our methodology.

## 1. Introduction

This paper is about semiparametric regression models with repeated measures when the primary goal is to estimate a population quantity such as mean, variance, probability, etc. We will construct estimators of these quantities which utilize the underlying semiparametric structure of the model and derive their limiting distribution.

The work is motivated by the following example: the Kenya hemoglobin data. The goal is to study the changes of hemoglobin over time during the first year of birth. The data set consists of 68 families with 2 children per family. For each child, 4 repeated measures are taken over time in the first year since birth: the time of visit varied from child to child. The factors include mother's age at child birth, child sex and placental parasitemia density (PDEN), a marker of malaria that could affect hemoglobin. To model these data, Lin and Carroll [2] considered a semiparametric model where the mother's age effect is modeled nonparametrically and (sex, PDEN) is modeled parametrically. The model is given by the repeated measures partially linear model

$$(1.1) \quad Y_{ijk} = X_{ijk}^{\mathrm{T}}\beta_0 + \theta_0(Z_{ij}) + \epsilon_{ijk},$$

where $i = 1, \ldots, n$, $j = 1, \ldots, m$ and $k = 1, \ldots, R$. Specifically, Lin and Carroll [2] set $Y_{ijk}$ = hemoglobin level of the $j^{\mathrm{th}}$ child in the $i^{\mathrm{th}}$ family at the $k^{\mathrm{th}}$ visit, $Z_{ij}$ =

---

[*]Supported by grants from the National Cancer Institute (CA57030, CA104620) and by the Texas A&M Center for Environmental and Rural Health via a grant from the National Institute of Environmental Health Sciences (P30-ES09106).
[1]Department of Statistics, Texas A&M University, College Station, TX 77843-3143, USA, e-mail: amaity@stat.tamu.edu; carroll@stat.tamu.edu
[2]School of Operations Research and Industrial Engineering Cornell University, Ithaca NY 14853, USA, e-mail: tva2@cornell.edu

mother's age at the $j^{\text{th}}$ child birth in the $i^{\text{th}}$ family, $X_{ijk} = \{\text{sex, logpden, month,} (\text{month}-4)_+\}$, where sex = 1 if female and 0 if male, logpden = log(PDEN+1), month denotes the age of the child at the $k^{\text{th}}$ visit and the function $f_+ = f$ if $f > 0$ and 0 if $f \leq 0$. As noted by Lin and Carroll [2], the covariates {month, $(\text{month}-4)_+\}$ are used to model the time effect as a piecewise linear function with a knot at 4 months: this trend being observed by preliminary analysis of the data set. Also assume that, $\epsilon_i = (\epsilon_{i11}, \epsilon_{i12}, \ldots, \epsilon_{imR})$ has a Normal$(0, \Sigma)$ distribution with $\Sigma = \sigma^2 I_{mR} + \rho\sigma^2(J_{mR} - I_{mR})$, where $I_{mR}$ is the $mR \times mR$ identity matrix and $J_{mR}$ is a $mR \times mR$ matrix with 1 as all elements. In this context, apart from the model components, namely $\beta_0$ and $\theta_0(\bullet)$, one is interested in other population level quantities such as mean hemoglobin, $E(Y)$ or various probabilities such as the proportion of six month old children who have hemoglobin measure above a given constant $c$, $\text{pr}(Y > c|\text{month} = 6)$, etc. Note that each of the above-mentioned population-level quantities can be written as functions of $X$, $Z$ and the model parameters. For example, because $Z$, sex and logpden are independent of month, we can write

$$(1.2) \quad \text{pr}(Y > c|\text{month} = a) = E\left(\Phi[\{X^{\text{T}}\beta_0 + \theta_0(Z) - c\}/\sigma]|\text{month} = a\right),$$

where $c$ and $a$ are given constants and $\Phi(\bullet)$ is the standard normal cdf. In general, we can consider any functional $\kappa_0 = E\{\mathcal{F}(\bullet)\}$ for some function $\mathcal{F}(\bullet)$: this of course includes such quantities as population mean, probabilities, etc.

The Kenya hemoglobin study example is a special case of a much more general framework, one we call semiparametric repeated measures regression modeling. Let $\widetilde{Z}_i = (Z_{i1}, \ldots, Z_{im})$ and let $\widetilde{Y}_i$ and $\widetilde{X}_i$ be the ensemble of responses and other covariates, respectively. Consider a semiparametric problem in which the loglikelihood function conditional on $(\widetilde{X}, Z_1, \ldots, Z_m)$ is $\mathcal{L}\{\widetilde{Y}, \widetilde{X}, \theta_0(Z_1), \ldots, \theta_0(Z_m), \mathcal{B}_0\}$, where $\mathcal{B}_0$ is the set of model parameters. Also as an important generalization, allow $\widetilde{Y}$ to be partially missing and let $\delta = 1$ if $\widetilde{Y}$ is observed and 0 otherwise. Suppose further that $\widetilde{Y}$ is missing at random, so that $\text{pr}(\delta = 1|\widetilde{Y}, \widetilde{X}, \widetilde{Z}) = \text{pr}(\delta = 1|\widetilde{X}, \widetilde{Z})$. If we define $\mathcal{L}_{\mathcal{B}}(\bullet)$ and $\mathcal{L}_{j\theta}(\bullet)$ to be the derivatives of the loglikelihood with respect to $\mathcal{B}$ and $\theta(Z_j)$, we have the properties that for $j = 1, 2, \ldots, m$,

$$E[\delta\mathcal{L}_{\mathcal{B}}\{\widetilde{Y}, \widetilde{X}, \theta_0(Z_1), \ldots, \theta_0(Z_m), \mathcal{B}_0\}|\widetilde{X}, Z_1 \ldots, Z_m] = 0;$$
$$E[\delta\mathcal{L}_{j\theta}\{\widetilde{Y}, \widetilde{X}, \theta_0(Z_1), \ldots, \theta_0(Z_m), \mathcal{B}_0\}|\widetilde{X}, Z_1 \ldots, Z_m] = 0.$$

The main goal of this paper is to estimate a population level quantity, namely

$$(1.3) \quad \kappa_0 = E[\mathcal{F}\{\widetilde{X}, \widetilde{Z}, \theta_0(Z_1), \ldots, \theta_0(Z_m), \mathcal{B}_0\}]$$

for an arbitrary but smooth function $\mathcal{F}(\bullet)$. For example, in (1.2) we have $\mathcal{F}(\bullet) = \Phi[\{X^{\text{T}}\beta_0 + \theta_0(Z) - c\}/\sigma]$ where last two components of $X$ are fixed as $a$ and $(a-4)_+$.

Although there are various papers about application of semiparametric methods in the context of longitudinal data (see for example Zeger and Diggle [5] and Zhang et al. [6]), to the best of our knowledge there is no literature on estimation of arbitrary population-level quantities in the repeated-measures context. See Maity et al. [3] for the case of no repeated measurement. However, estimation of $\mathcal{B}_0$ and $\theta_0(\bullet)$ is described in Lin and Carroll [2]. If we call their estimates $\widehat{\mathcal{B}}$ and $\widehat{\theta}(\bullet, \widehat{\mathcal{B}})$, then we propose to estimate the population-level quantity as

$$(1.4) \quad \widehat{\kappa}_{\text{semi}} = n^{-1}\sum_{i=1}^{n}\mathcal{F}\{\widetilde{X}_i, \widetilde{Z}_i, \widehat{\theta}(Z_{i1}, \widehat{\mathcal{B}}), \ldots, \widehat{\theta}(Z_{im}, \widehat{\mathcal{B}}), \widehat{\mathcal{B}}\}.$$

In this paper, we are interested in the asymptotic behavior of $\widehat{\kappa}_{\text{semi}}$.

This paper is organized as follows. In Section 2, we discuss the general semiparametric problem with loglikelihood $\mathcal{L}(\bullet)$ and a general goal of estimating $\kappa_0$. We also derive the asymptotic distribution of (1.4). In Section 3, we demonstrate our method by applying it to analyze the Kenya hemoglobin study data. All technical results are given in an Appendix.

## 2. Semiparametric models with multiple components

### 2.1. Model component estimation

Estimation of $\mathcal{B}_0$ and $\theta_0(\bullet)$ using profile likelihood and backfitting methods is discussed in detail in Lin and Carroll [2] when there are no missing data. With missing data, the method works as follows. Let $K(\bullet)$ be a smooth symmetric density function with bounded support and variance 1.0, let $h$ be a bandwidth and let $K_h(z) = h^{-1}K(z/h)$. Define $G_{ij}(z,h) = \{1, (Z_{ij} - z)/h\}$. For any fixed $\mathcal{B} = \mathcal{B}_*$, let $\widehat{\theta}_c(\bullet)$ be the current estimate in an iterative procedure and $(\widehat{\alpha}_0, \widehat{\alpha}_1)$ be the solution of the following estimating equation:

$$0 = \sum_{i=1}^{n}\sum_{j=1}^{m} \delta_i K_h(Z_{ij} - z) G_{ij}(z,h)$$
$$\times \mathcal{L}_{j\theta}\left\{\widetilde{Y}_i, \widetilde{X}_i, \widehat{\theta}_c(Z_{i1}, \mathcal{B}_*), \ldots, \alpha_0 + \alpha_1(Z_{ij} - z), \ldots, \widehat{\theta}_c(Z_{im}, \mathcal{B}_*)\right\}.$$

We then update $\widehat{\theta}(\bullet)$ as $\widehat{\theta}(z, \mathcal{B}_*) = \widehat{\alpha}_0$.

To estimate $\mathcal{B}_0$, a backfitting algorithm can be used. In the iterated backfitting algorithm, suppose the current estimate is $\mathcal{B}_*$. Then the backfitting estimator of $\mathcal{B}_0$ modified for missing responses is updated by maximizing in $\mathcal{B}$ the function

$$\sum_{i=1}^{n} \delta_i \mathcal{L}\left\{\widetilde{Y}_i, \widetilde{X}_i, \widehat{\theta}(Z_{i1}, \mathcal{B}_*), \ldots, \widehat{\theta}(Z_{im}, \mathcal{B}_*), \mathcal{B}\right\},$$

i.e., solving for $\mathcal{B}$

$$0 = \sum_{i=1}^{n} \delta_i \mathcal{L}_{\mathcal{B}}\left\{\widetilde{Y}_i, \widetilde{X}_i, \widehat{\theta}(Z_{i1}, \mathcal{B}_*), \ldots, \widehat{\theta}(Z_{im}, \mathcal{B}_*), \mathcal{B}\right\}.$$

The final estimates are obtained by iterating the process until convergence.

### 2.2. Estimation of general population-level summaries

We take advantage of the results about the asymptotic expansions for $\widehat{\mathcal{B}}$ and $\widehat{\theta}(\bullet)$ provided in Lin and Carroll [2], with the modification of incorporating the missing data indicators. Let $\mathcal{L}_{jk\theta}(\bullet) = \partial \mathcal{L}_{j\theta}\{Y, X, \theta_0(Z_1), \ldots, \theta_0(Z_m), \mathcal{B}_0\}/\partial \theta(Z_k)$. Make the definitions $\Omega(z) = \sum_{j=1}^{m} f_j(z) E\{\delta \mathcal{L}_{jj\theta}(\bullet)|Z_j = z\}$ and

$$\mathcal{A}(B, z_1, z_2) = \sum_{j=1}^{m}\sum_{k \neq j}^{m} f_j(z_1) E\left\{\delta \mathcal{L}_{jk\theta}(\bullet) B(Z_k, z_2)/\Omega(Z_k)|Z_j = z_1\right\};$$

$$Q(z_1, z_2) = \sum_{j=1}^{m}\sum_{k \neq j}^{m} f_{jk}(z_1, z_2) E\left\{\delta \mathcal{L}_{jk\theta}(\bullet)|Z_j = z_1, Z_k = z_2\right\}/\Omega(z_2),$$

where $f_j(z)$ is the density of $Z_j$ and $f_{jk}(z_1, z_2)$ is the bivariate density of $(Z_j, Z_k)$, assumed to have bounded support and are positive on the support. Let $\mathcal{G}(z_1, z_2)$ be the solution to

$$\mathcal{G}(z_1, z_2) = \mathcal{Q}(z_1, z_2) - \mathcal{A}(\mathcal{G}, z_1, z_2).$$

Define $\mathcal{L}_{j\theta\mathcal{B}}(\bullet) = \partial \mathcal{L}_{j\theta}\{\widetilde{Y}, \widetilde{X}, \theta_0(Z_1), \ldots, \theta_0(Z_m), \mathcal{B}_0\}/\partial \mathcal{B}$ and

$$
\begin{aligned}
\epsilon_{ij}^{\#}(\theta, \mathcal{B}) &= \mathcal{L}_{j\theta\mathcal{B}}\{\widetilde{Y}_i, \widetilde{X}_i, \theta(Z_{i1}), \ldots, \theta(Z_{im}), \mathcal{B}\} \\
&\quad + \sum_{k=1}^{m} \mathcal{L}_{jk\theta}\{\widetilde{Y}_i, \widetilde{X}_i, \theta(Z_{i1}), \ldots, \theta(Z_{im}), \mathcal{B}\}\theta_{\mathcal{B}}(Z_{ik}, \mathcal{B}),
\end{aligned}
\tag{2.1}
$$

where $\theta_{\mathcal{B}}(z, \mathcal{B}_0)$ satisfies

$$0 = \sum_{j=1}^{m} f_j(z) E\{\delta \epsilon_{ij}^{\#}(\theta_0, \mathcal{B}_0) | Z_j = z\}. \tag{2.2}$$

Define

$$\epsilon_i = \mathcal{L}_{i\mathcal{B}}(\bullet) + \sum_{j=1}^{m} \mathcal{L}_{ij\theta}(\bullet)\theta_{\mathcal{B}}(Z_{ij}, \mathcal{B}_0);$$

$$\mathcal{M}_1 = E(\delta \epsilon \epsilon^{\mathrm{T}}) = -E\{\delta_i \mathcal{L}_{\mathcal{B}\mathcal{B}}(\bullet) + \sum_{j=1}^{m} \delta_i \mathcal{L}_{j\theta\mathcal{B}}(\bullet)\theta_{\mathcal{B}}^{\mathrm{T}}(Z_j, \mathcal{B}_0)\}.$$

Define the derivatives of $\mathcal{F}(\bullet)$ in a similar way to $\mathcal{L}(\bullet)$ and make the following definitions:

$$\mathcal{M}_2 = E[\mathcal{F}_{\mathcal{B}}(\bullet) + \sum_{j=1}^{m} \mathcal{F}_{j\theta}(\bullet)\theta_{\mathcal{B}}(Z_j, \mathcal{B}_0)];$$

$$C_1(z) = -\sum_{j=1}^{m} f_j(z) E\{\mathcal{F}_{j\theta}(\bullet) | Z_j = z\}/\Omega(z);$$

$$C_2(z) = E[\sum_{j=1}^{m} E\{\mathcal{F}_{j\theta}(\bullet) | Z_j\}\mathcal{G}(Z_j, z)/\Omega(Z_j)];$$

$$D(\bullet) = \sum_{j=1}^{m} \mathcal{L}_{j\theta}(\bullet)\{C_1(Z_j) + C_2(Z_j)\}.$$

Then we have the following result.

**Theorem 2.1.** *Assume that* $(\widetilde{Y}_i, \delta_i, \widetilde{X}_i, Z_{i1}, \ldots, Z_{im}), i = 1, 2, \ldots, n$ *are independently and identically distributed and* $nh^4 \to 0$. *Then, to terms of order* $o_p(1)$,

$$n^{1/2}(\widehat{\kappa}_{\mathrm{semi}} - \kappa_0) \approx n^{-1/2} \sum_{i=1}^{n} \{\mathcal{F}_i(\bullet) - \kappa_0 + \mathcal{M}_2^{\mathrm{T}} \mathcal{M}_1^{-1} \delta_i \epsilon_i + \delta_i D_i(\bullet)\} \tag{2.3}$$

$$\to \text{Normal } (0, \mathcal{V}_\kappa), \tag{2.4}$$

where $\mathcal{V}_\kappa = E\{\mathcal{F}(\bullet) - \kappa_0\}^2 + \mathcal{M}_2^{\mathrm{T}} \mathcal{M}_1^{-1} \mathcal{M}_2 + E\{\delta D^2(\bullet)\}$.

Note that we have assumed that $\kappa_0 = E[\mathcal{F}\{\widetilde{X}, \widetilde{Z}, \theta_0(Z_1), \ldots, \theta_0(Z_m), \mathcal{B}_0\}]$ and $\mathcal{F}(\bullet)$ are scalar, which in principle excludes the estimation of quantities such as the population variance and standard deviation. However, the result can be readily generalized to handle this case. In general, suppose $\kappa_1 = E\{\mathcal{F}_1(\bullet)\}$ and $\kappa_2 = E\{\mathcal{F}_2(\bullet)\}$ for some smooth functions $\mathcal{F}_1(\bullet)$ and $\mathcal{F}_2(\bullet)$, and are estimated by $\widehat{\kappa}_1$ and $\widehat{\kappa}_2$, respectively. Suppose that we are interested in $\kappa_{\text{gen}} = g(\kappa_1, \kappa_2)$ for some function $g(\bullet, \bullet)$. Then we have the following result:

**Corollary 2.1.** *Let $\kappa_{\text{gen}}$ be estimated by $\widehat{\kappa}_{\text{gen}} = g(\widehat{\kappa}_1, \widehat{\kappa}_2)$ and that $g(\bullet, \bullet)$ is differentiable with respect to both of its arguments. Let $g_1(\bullet, \bullet)$ and $g_2(\bullet, \bullet)$ be the derivatives of $g$ with respect to its first and second argument, respectively. Then we have*

$$\begin{aligned}
n^{1/2}(\widehat{\kappa}_{\text{gen}} - \kappa_{\text{gen}}) &= n^{1/2}g_1(\kappa_1, \kappa_2)(\widehat{\kappa}_1 - \kappa_1) + n^{1/2}g_2(\kappa_1, \kappa_2)(\widehat{\kappa}_2 - \kappa_2) + o_p(1) \\
&\to Normal(0, \mathcal{V}_{\text{gen}}),
\end{aligned}$$

*with* $\mathcal{V}_{\text{gen}} = g_1^2(\kappa_1, \kappa_2)\mathcal{V}_1 + g_2^2(\kappa_1, \kappa_2)\mathcal{V}_2 + 2g_1(\kappa_1, \kappa_2)g_2(\kappa_1, \kappa_2)\mathcal{V}_{12}$, *where $\mathcal{V}_1$ and $\mathcal{V}_2$ are asymptotic variances of $\widehat{\kappa}_1$ and $\widehat{\kappa}_2$, respectively, computed using (2.4) and $\mathcal{V}_{12} = E[\{\mathcal{F}_1(\bullet) - \kappa_1\}\{\mathcal{F}_2(\bullet) - \kappa_2\}]$.*

**Remark 2.1.** We note that if $\kappa_0$ were a function of only $\mathcal{B}_0$, then a bandwidth proportional to $h \sim n^{-1/5}$ can be used, because in this case undersmoothing is not necessary. The reason for the undersmoothing in the general case is the inclusion of $\theta_0(\bullet)$ in $\kappa_0$. Using the results from Lin and Carroll [2] we have

$$\begin{aligned}
\widehat{\theta}(z) - \theta_0(z) &= O_p(h^2) - n^{-1}\sum_{i=1}^{n}\sum_{j=1}^{m}\delta_i K_h(Z_{ij} - z)\mathcal{L}_{ij\theta}(\bullet)/\Omega(z) \\
&\quad + n^{-1}\sum_{i=1}^{n}\sum_{j=1}^{m}\delta_i \mathcal{L}_{ij\theta}(\bullet)\mathcal{G}(z, Z_{ij})/\Omega(z) + o_p(n^{-1/2}).
\end{aligned}$$

The restriction $nh^4 \to 0$ removes the $O_p(h^2)$ bias term. It is suggested in Sepanski et al. [4] that in many semiparametric problems, the optimal bandwidth for estimating parameters such as $\mathcal{B}_0$ is of the order $n^{-1/3}$, which of course satisfies $nh^4 \to 0$. In the case of no repeated measures, i.e., when $m = 1$, Maity et al. [3] find that they have good experience numerically by first estimating the bandwidth via likelihood crossvalidation, which will be of order $n^{-1/5}$, and then multiplying it by $n^{-2/15}$ to get the bandwidth to be of order $n^{-1/3}$.

**Remark 2.2.** To estimate the asymptotic variance of $\kappa_0$, one can use the bootstrap method. While the justification of use of parametric bootstrap for estimation of model parameters is provided in Chen et al. [1], we conjecture that the bootstrap works for $\kappa_0$ as well. One can also use the plug-in method where one replaces each term in the variance expression in (2.4) by their consistent estimators. Constructing consistent estimators for the first two terms are fairly straightforward where one merely replaces all the expectations by sums in that expression and all the regression functions by kernel estimates. However, the main difficulty lies in the estimation of the third term where one has to solve the functional equation $\mathcal{G}(z_1, z_2) = Q(z_1, z_2) - \mathcal{A}(\mathcal{G}, z_1, z_2)$ for $\mathcal{G}(\bullet, \bullet)$, which can be numerically difficult.

### 2.3. General functions of the response

One interesting and important scenario is when $\kappa_0$ can be constructed using only the responses. Suppose that $\kappa_0 = E\{\mathcal{G}(\widetilde{Y})\}$. Also define $\mathcal{F}\{\widetilde{X}, \widetilde{Z}, \theta_0(Z_1), \ldots, \theta_0(Z_m),$ $\mathcal{B}_0\} = E\{\mathcal{G}(\widetilde{Y})|\widetilde{X},\widetilde{Z}\}$. In the presence of missing data, there are various estimators one can use. We will discuss two of these estimators which are based upon different constructions for estimating the missing data process.

Let $\widehat{\mathcal{F}}_i(\bullet) = \mathcal{F}\{\widetilde{X}_i,\widetilde{Z}_i,\widehat{\theta}(Z_{i1},\widehat{\mathcal{B}}),\ldots,\widehat{\theta}(Z_{im},\widehat{\mathcal{B}}),\widehat{\mathcal{B}}\}$. To form the first estimator, we impute $\mathcal{G}(\widetilde{Y}_i)$ by $\widehat{\mathcal{F}}_i(\bullet)$ for every missing $\widetilde{Y}_i$. We estimate $\kappa_0$ by

$$\widehat{\kappa}_1 = n^{-1}\sum_{i=1}^n \{\delta_i \mathcal{G}(\widetilde{Y}_i) + (1-\delta_i)\widehat{\mathcal{F}}_i(\bullet)\}.$$

Note that in absence of missing data, $\widehat{\kappa}_1$ does not use the semiparametric model and is consistent. If the responses are missing at random, then $\widehat{\kappa}_1$ is consistent if the semiparametric model specification is correct.

The second estimator depends on the underlying structure of the missing data process where we assume a parametric formulation for estimating $\mathrm{pr}(\delta = 1|\widetilde{Y},\widetilde{X},\widetilde{Z}) = \pi(\widetilde{X},\widetilde{Z},\zeta)$, where $\zeta$ is an unknown parameter estimated by standard logistic regression of $\delta$ on $(\widetilde{X},\widetilde{Z})$. Then construct the estimator

$$\widehat{\kappa}_2 = n^{-1}\sum_{i=1}^n \left[\frac{\delta_i}{\pi(\widetilde{X}_i,\widetilde{Z}_i,\widehat{\zeta})}\mathcal{G}(\widetilde{Y}_i) + \left\{1 - \frac{\delta_i}{\pi(\widetilde{X}_i,\widetilde{Z}_i,\widehat{\zeta})}\right\}\widehat{\mathcal{F}}_i(\bullet)\right].$$

This estimator has the double-robustness property that if either the semiparametric model specification for $(\theta_0,\mathcal{B}_0)$ or the parametric model for $\pi(\widetilde{X},\widetilde{Z},\zeta)$ is correct, $\widehat{\kappa}_2$ is consistent and has an asymptotic normal distribution.

If both the models are correct then the following results are true. The proofs are given in the Appendix.

**Lemma 2.1.** *Define*

$$\mathcal{M}_{2,\mathrm{imp}} = E[(1-\delta)\{\mathcal{F}_\mathcal{B}(\bullet) + \sum_{j=1}^m \mathcal{F}_{j\theta}(\bullet)\theta_\mathcal{B}(Z_j,\mathcal{B}_0)\}];$$

$$C_{1,\mathrm{imp}}(z) = -\sum_{j=1}^m f_j(z)E\{(1-\delta)\mathcal{F}_{j\theta}(\bullet)|Z_j=z\}/\Omega(z);$$

$$C_{2,\mathrm{imp}}(z) = E[\sum_{j=1}^m E\{(1-\delta)\mathcal{F}_{j\theta}(\bullet)|Z_j\}\mathcal{G}(Z_j,z)/\Omega(Z_j)];$$

$$D_{\mathrm{imp}}(\bullet) = \sum_{j=1}^m \mathcal{L}_{j\theta}(\bullet)\{C_{1,\mathrm{imp}}(Z_j) + C_{2,\mathrm{imp}}(Z_j)\}.$$

*Assume that $nh^4 \to 0$. Then, to terms of order $o_p(n^{-1/2})$,*

$$\widehat{\kappa}_1 - \kappa_0 \approx n^{-1}\sum_{i=1}^n \Big\{\delta_i\mathcal{G}(\widetilde{Y}_i) + (1-\delta_i)\mathcal{F}_i(\bullet) - \kappa_0$$

(2.5)
$$+ \mathcal{M}_{2,\mathrm{imp}}^{\mathrm{T}}\mathcal{M}_1^{-1}\delta_i\epsilon_i + \delta_i D_{i,\mathrm{imp}}(\bullet)\Big\}.$$

**Lemma 2.2.** *Define $\pi_\zeta(\widetilde{X},\widetilde{Z},\zeta) = \partial\pi(\widetilde{X},\widetilde{Z},\zeta)/\partial\zeta$. Assume that $n^{1/2}(\widehat{\zeta}-\zeta) = n^{-1/2}\sum_{i=1}^n \psi_{i\zeta}(\bullet) + o_p(n^{-1/2})$ with $E\{\psi_\zeta(\bullet)|\widetilde{X},\widetilde{Z}\} = 0$. Also assume that $nh^4 \to 0$.*

*Then, to terms of order $o_p(n^{-1/2})$,*

$$\widehat{\kappa}_2 - \kappa_0 \approx n^{-1} \sum_{i=1}^{n} \left[ \frac{\delta_i}{\pi(\widetilde{X}_i, \widetilde{Z}_i, \zeta)} \{\mathcal{G}(\widetilde{Y}_i) - \kappa_0\} \right.$$
$$\left. + \left\{1 - \frac{\delta_i}{\pi(\widetilde{X}_i, \widetilde{Z}_i, \zeta)}\right\} \{\mathcal{F}_i(\bullet) - \kappa_0\} \right]. \tag{2.6}$$

**Remark 2.3.** The expansions (2.5) and (2.6) show that $\widehat{\kappa}_1$ and $\widehat{\kappa}_2$ are asymptotically normally distributed. One can show that the asymptotic variances are given as

$$\mathcal{V}_{\kappa_1} = \mathrm{var}\left\{\delta\mathcal{G}(\widetilde{Y}) + (1-\delta)\mathcal{F}(\bullet) + \mathcal{M}_{2,\mathrm{imp}}^{\mathrm{T}}\mathcal{M}_1^{-1}\delta\epsilon + \delta D_{\mathrm{imp}}(\bullet)\right\};$$
$$\mathcal{V}_{\kappa_2} = \mathrm{var}\left[\frac{\delta}{\pi(\widetilde{X},\widetilde{Z},\zeta)}\mathcal{G}(\widetilde{Y}) + \left\{1 - \frac{\delta}{\pi(\widetilde{X},\widetilde{Z},\zeta)}\right\}\mathcal{F}(\bullet)\right],$$

respectively, from which estimates are readily derived. Note that the estimation of $\mathcal{B}_0$ and $\theta_0(\bullet)$ has no impact on the asymptotic distribution of $\widehat{\kappa}_2$. Hence it has the distinct advantage that its limiting distribution does not involve an integral equation.

## 3. Application

### 3.1. Motivating example

We applied our method to analyze the Kenya hemoglobin data which is described in Section 1. The model is given by (1.1) where $n = 68$, $m = 2$ and $R = 4$. Comparing this to the general model, we see that $\widetilde{Y}_i = (y_{i11}, \ldots, y_{i24})^{\mathrm{T}}$, $\widetilde{X}_i = (X_{i11}, \ldots, X_{i24})^{\mathrm{T}}$ and $Z_{i1} = Z_{i2} = Z_{i3} = Z_{i4}$ and $Z_{i5} = Z_{i6} = Z_{i7} = Z_{i8}$. Also, $\mathcal{B} = (\beta^{\mathrm{T}}, \sigma)^{\mathrm{T}}$. Let $e_4$ be a vector of ones of length 4. The loglikelihood function, conditional of $\widetilde{X}$ and $(Z_1, Z_2)$, is that of a multivariate Gaussian distribution with mean $\widetilde{X}^{\mathrm{T}}\beta_0 + \{\theta_0(Z_1)e_4^{\mathrm{T}}, \theta_0(Z_2)e_4^{\mathrm{T}}\}^{\mathrm{T}}$ and covariance matrix $\Sigma$.

Recall that we are interested in the proportion of all children of a given age $a$ who have their hemoglobin measure above any given constant $c$, which is defined as

$$\kappa_0 = E\Big(\Phi[\{X^{\mathrm{T}}\beta_0 + \theta_0(Z) - c\}/\sigma]|\mathrm{month} = a\Big). \tag{3.1}$$

Note that the joint density of $(Z, \mathrm{sex}, \mathrm{logpden})$ is independent of "month" and hence the conditional expectation in (3.1) can be written as $E(\Phi[\{X^{\mathrm{T}}\beta_0 + \theta_0(Z) - c\}/\sigma])$ with the last two components of $X$ being fixed as $a$ and $(a-4)_+$. Using this fact, we estimate $\kappa_0$ by

$$\widehat{\kappa}_{\mathrm{semi}} = \{nmR\}^{-1} \sum_{i=1}^{n}\sum_{j=1}^{m}\sum_{k=1}^{R} \Phi[\{X_{ijk}^{\mathrm{T}}\widehat{\beta} + \widehat{\theta}(Z_{ij}, \widehat{\beta}) - c\}/\widehat{\sigma}], \tag{3.2}$$

where the last two components of $X_{ijk}$ are fixed as $a$ and $(a-4)_+$. Also note that there are no missing responses involved in this setup.

We proceed to estimate $\mathcal{B}_0$ and $\theta_0(\bullet)$ as indicated in Lin and Carroll [2] and use these estimates to compute (3.2). We estimate $\widehat{\kappa}_{\mathrm{semi}}$ for $a = 3, 6$ and for different values of $c$. We estimate the limiting variances via the bootstrap and use those to construct 95% confidence intervals. The results are displayed in Figure 1. It is interesting to note that the two estimates differ significantly and a decrease in proportion is evident as the age of the children increases. We can observe form Figure 1 that for a given hemoglobin level, the proportion of children having a higher hemoglobin level at age three is significantly less than at age six, indicating an increase in risk of getting anemia (low hemoglobin level) as children grow older. To observe this effect, we set the hemoglobin level, $c = 10$ and estimated the survival function over a grid of value for age, $a$. Since the time effect is modeled as a piecewise linear function with a knot at 4 months, we plot the estimated survival function in two different plots in Figure 2: left panel displays the survival function as age varies in $[0.1, 4]$ and right panel plots the survival function when age varies in $(4, 11]$. The decrease in hemoglobin level is much faster for children below age 4 months than children older than 4 months, after which the survival function becomes flatter.

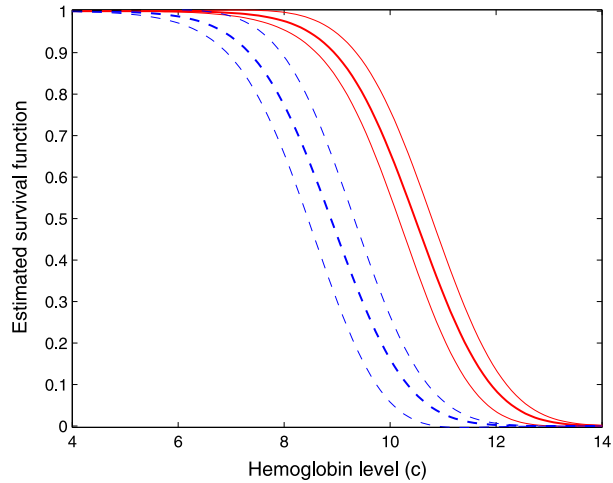

Fig 1. *Results for the Kenya hemoglobin data example. Plotted are estimated survival function (1 − the cdf) of hemoglobin measure of children of a given age a, along with 95% confidence intervals. Dashed lines correspond to six-month-old children and solid lines correspond to three-month-old children.*

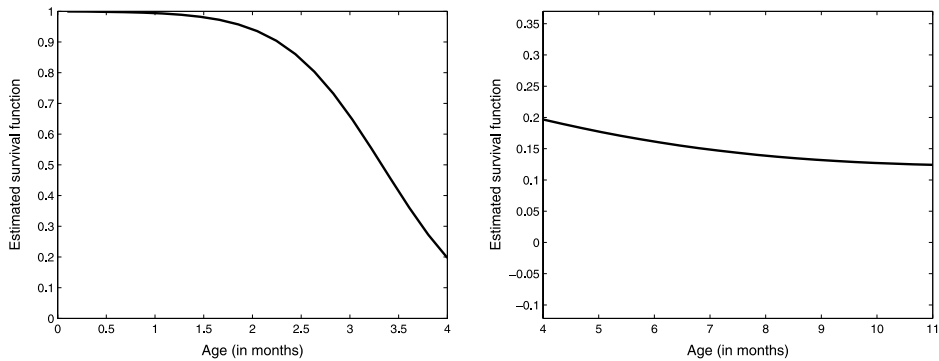

Fig 2. *Results for the Kenya hemoglobin data example. Plotted are estimated survival function (1 − the cdf) for a fixed hemoglobin level c = 10 as age varies from 0.1 months to 4 months (left panel) and 4 months to 11 months (right panel).*

Note that one can construct other estimators for $\kappa_0$, for example, the sample average of hemoglobin measures of only those children who have age $a$. However, it is important to note that the naive estimator does not use the information that "month" is independent of $(Z, \text{sex}, \text{logpden})$. On the other hand, this vital information is used to construct the semiparametric estimator thus allowing $\widehat{\kappa}$ to be more efficient than the naive estimator. Also, note that in this example the use of the usual naive estimators may not be justified because the number of children of any given age $a$ may be very small. In fact, there are only six children with age $a = 3$ and only one child with age $a = 6$ in the data set we used. This problem does not arise when $\widehat{\kappa}_{\text{semi}}$ is used because of the fact that it uses all the data to estimate $\kappa_0$ rather than using only those children who have age $a$.

### 3.2. Simulation study

We conducted a small simulation study to observe the performance of our method. We considered the model given in (1.1) with $n = 100$ clusters with six measurements per cluster. Specifically, we set $m = 2$ and $R = 3$ with $Z_{i1} = Z_{i2} = Z_{i3}$ and $Z_{i4} = Z_{i5} = Z_{i6}$. We assume that $\epsilon_i$ has a Normal$(0, \Sigma)$ distribution with $\Sigma = \sigma^2 I_{mR} + \rho \sigma^2 (J_{mR} - I_{mR})$. We set $\sigma^2 = 1$ and $\rho = 0.4$. We generated $Z$ from a uniform$(0, 1)$ distribution with the true function $\theta_0(z) = \sin(8z - 1)$. The individual components of $X$ are generated from uniform$(0, 1)$ distribution, independent of each other and $Z$. The true value of $\beta_0 = [1, 1]^{\mathrm{T}}$.

In this setup, our target was to estimate

$$\begin{aligned} \kappa_0(c) &= \text{pr}(Y > c | X_1 = 0.5) \\ &= E\Big(\Phi[\{X^{\mathrm{T}}\beta_0 + \theta_0(Z) - c\}/\sigma] | X_1 = 0.5\Big) \end{aligned}$$

for a given value of $c$, where $X_1$ denotes the first component of $X$. We estimate $\kappa_0(c)$ by

$$\widehat{\kappa}(c) = \{nmR\}^{-1} \sum_{i=1}^{n} \sum_{j=1}^{m} \sum_{k=1}^{R} \Phi[\{X_{ijk}^{\mathrm{T}}\widehat{\beta} + \widehat{\theta}(Z_{ij}, \widehat{\beta}) - c\}/\widehat{\sigma}],$$

where the first component of $X_{ijk}$ is fixed to be 0.5. Note again that in definition of $\widehat{\kappa}(c)$, we used independence between $X$ and $Z$ to our advantage.

We employed the method in Lin and Carroll [2] with the Epanechnikov kernel function to obtain $\widehat{\beta}$ and $\widehat{\theta}(\bullet)$. Using these estimates we calculated $\widehat{\kappa}(c)$ over a grid of values for $c$. Also the true survival function was computed analytically. We generated 1000 data sets, with results displayed in Figure 3. It is evident that the estimated survival function given by our method is very close to the true function, as one would expect for large sample size. Note that $X$ is generated from a continuous distribution and hence the use of naive estimator is not justified in this case.

### 4. Discussion

In this paper, we considered a general class of semiparametric models, with responses missing at random, where the primary goal was to estimate population-level quantities $\kappa_0$ such as the mean, probabilities, etc. For general semiparametric regression models, the asymptotic distribution of a plug-in estimator of $\kappa_0$ was derived.

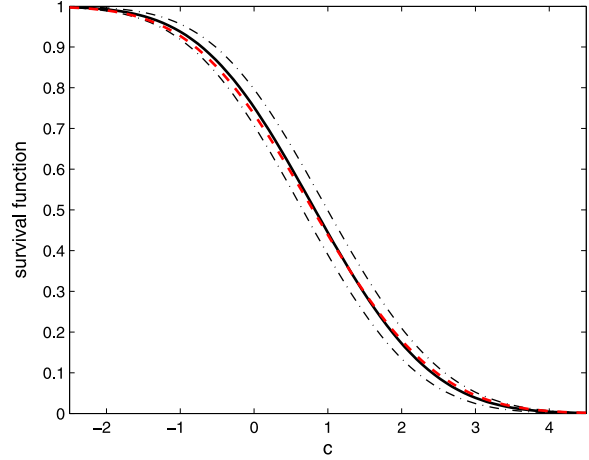

FIG 3. *Results for the simulation study based on 1000 simulated data sets. Plotted are the true survival function (dashed line), the mean estimated survival function (1 − the cdf) (solid line), along with a 95% confidence range(dash-dotted line).*

We have considered the case that the unknown function $\theta_0(Z)$ was a scalar function of a scalar argument. The results though readily extend to the case of a multivariate function of a scalar argument.

## Appendix A: Sketch of technical arguments

### A.1. Sketch of Theorem 2.1

We first show (2.3). First note that $\mathcal{L}$ is a loglikelihood function conditioned on $(X, Z)$, so that we have

$$
\begin{aligned}
E\{\delta \mathcal{L}_{jk\theta}(\bullet)|\widetilde{X}, \widetilde{Z}\} &= -E\{\delta \mathcal{L}_{j\theta}(\bullet)\mathcal{L}_{k\theta}(\bullet)|\widetilde{X}, \widetilde{Z}\}; \\
E\{\delta \mathcal{L}_{j\theta\mathcal{B}}(\bullet)|\widetilde{X}, \widetilde{Z}\} &= -E\{\delta \mathcal{L}_{j\theta}(\bullet)\mathcal{L}_{\mathcal{B}}(\bullet)|\widetilde{X}, \widetilde{Z}\}.
\end{aligned}
\tag{A.1}
$$

Also, it follows from Lin and Carroll [2], and $h = o_p(n^{-1/4})$ that

$$
\begin{aligned}
\widehat{\theta}(z, \mathcal{B}_0) - \theta_0(z) &= -n^{-1}\sum_{i=1}^{n}\sum_{j=1}^{m}\delta_i K_h(Z_{ij} - z)\mathcal{L}_{ij\theta}(\bullet)/\Omega(z) \\
&\quad + n^{-1}\sum_{i=1}^{n}\sum_{j=1}^{m}\delta_i \mathcal{L}_{ij\theta}(\bullet)\mathcal{G}(z, Z_{ij})/\Omega(z) + o_p(n^{-1/2});
\end{aligned}
\tag{A.2}
$$

$$
\widehat{\mathcal{B}} - \mathcal{B}_0 = \mathcal{M}_1^{-1} n^{-1}\sum_{i=1}^{n}\delta_i \epsilon_i + o_p(n^{-1/2}).
\tag{A.3}
$$

By a Taylor expansion,

$$
\begin{aligned}
n^{1/2}(\widehat{\kappa} - \kappa_0) &= n^{-1/2} \sum_{i=1}^{n} \bigg[ \mathcal{F}_i(\bullet) - \kappa_0 \\
&\quad + \{\mathcal{F}_{i\mathcal{B}}(\bullet) + \sum_{j=1}^{m} \mathcal{F}_{ij\theta}(\bullet)\theta_{\mathcal{B}}(Z_{ij}, \mathcal{B}_0)\}^{\mathrm{T}} (\widehat{\mathcal{B}} - \mathcal{B}_0) \\
&\quad + \sum_{j=1}^{m} \mathcal{F}_{ij\theta}(\bullet)\{\widehat{\theta}(Z_{ij}, \mathcal{B}_0) - \theta_0(Z_{ij})\}\bigg] + o_p(1) \\
&= n^{-1/2} \sum_{i=1}^{n} \{\mathcal{F}_i(\bullet) - \kappa_0\} + \mathcal{M}_2^{\mathrm{T}} n^{1/2} (\widehat{\mathcal{B}} - \mathcal{B}_0) \\
&\quad + n^{-1/2} \sum_{i=1}^{n} \sum_{j=1}^{m} \mathcal{F}_{ij\theta}(\bullet) \{\widehat{\theta}(Z_{ij}, \mathcal{B}_0) - \theta_0(Z_{ij})\} + o_p(1).
\end{aligned}
$$

Because $nh^4 \to 0$, using (A.2), we see

$$
\begin{aligned}
&n^{-1/2} \sum_{i=1}^{n} \sum_{j=1}^{m} \mathcal{F}_{ij\theta}(\bullet) \{\widehat{\theta}(Z_{ij}, \mathcal{B}_0) - \theta_0(Z_{ij})\} \\
&= -n^{-1/2} \sum_{i=1}^{n} \sum_{j=1}^{m} \mathcal{F}_{ij\theta}(\bullet) [n^{-1} \sum_{k=1}^{n} \sum_{\ell=1}^{m} \delta_k K_h(Z_{k\ell} - Z_{ij}) \mathcal{L}_{k\ell\theta}(\bullet) / \Omega(Z_{ij})] \\
&\quad + n^{-1/2} \sum_{i=1}^{n} \sum_{j=1}^{m} \mathcal{F}_{ij\theta}(\bullet) [n^{-1} \sum_{k=1}^{n} \sum_{\ell=1}^{m} \delta_k \mathcal{L}_{k\ell\theta}(\bullet) \mathcal{G}(Z_{ij}, Z_{k\ell}) / \Omega(Z_{ij})] + o_p(1) \\
&= n^{-1/2} \sum_{k=1}^{n} \sum_{\ell=1}^{m} \delta_k \mathcal{L}_{k\ell\theta}(\bullet) C_1(Z_{k\ell}) \\
&\quad + n^{-1/2} \sum_{k=1}^{n} \sum_{\ell=1}^{m} \delta_k \mathcal{L}_{k\ell\theta}(\bullet) C_2(Z_{k\ell}) + o_p(1) \\
&= n^{-1/2} \sum_{i=1}^{n} \delta_i D_i(\bullet) + o_p(1).
\end{aligned}
$$

Result (2.3) now follows from (A.3).

Next, we show (2.4). Recall that, $\epsilon = \mathcal{L}_{\mathcal{B}}(\bullet) + \sum_{j=1}^{m} \mathcal{L}_{j\theta}(\bullet) \theta_{\mathcal{B}}(Z_j, \mathcal{B}_0)$ and $D(\bullet) = \sum_{j=1}^{m} \mathcal{L}_{j\theta}(\bullet) \{C_1(Z_j) + C_2(Z_j)\}$. We use (2.1), (2.2) and (A.1) to derive that

$$
\begin{aligned}
E\{\delta \epsilon D(\bullet)\} &= -\sum_{j=1}^{m} E\bigg[\bigg\{\delta \mathcal{L}_{j\theta\mathcal{B}}(\bullet) + \sum_{k=1}^{m} \delta \mathcal{L}_{jk\theta}(\bullet) \theta_{\mathcal{B}}(Z_k, \mathcal{B}_0)\bigg\} \\
&\quad \times \{C_1(Z_j) + C_2(Z_j)\}\bigg] \\
&= -E\bigg[\sum_{j=1}^{m} E\{\delta \epsilon_j^{\#}(\theta_0, \mathcal{B}_0) | Z_j\} \{C_1(Z_j) + C_2(Z_j)\}\bigg] \\
&= 0.
\end{aligned}
$$

Also using the facts that $E\{D(\bullet)|\widetilde{X}, \widetilde{Z}\} = 0$ and $E(\epsilon|\widetilde{X}, \widetilde{Z}) = 0$, it is readily seen

that

$$E[\{\mathcal{F}(\bullet) - \kappa_0\}D(\bullet)] = 0;$$
$$E[\{\mathcal{F}(\bullet) - \kappa_0\}\epsilon] = 0,$$

which in turn prove that all three terms in (2.3) are uncorrelated. Then (2.4) follows easily from the central limit theorem.

### A.2. Sketch of Corollary 2.1

By a Taylor's series expansion,

$$\widehat{\kappa}_{\text{gen}} - \kappa_{\text{gen}} = \{g_1(\kappa_1, \kappa_2)(\widehat{\kappa}_1 - \kappa_1) + g_2(\kappa_1, \kappa_2)(\widehat{\kappa}_2 - \kappa_2)\} + o_p(n^{-1/2}).$$

Then, the normality follows from central limit theorem. The variance calculation is straight forward since (2.3) together with the fact $E(D|X, Z) = E(\epsilon|X, Z) = 0$ implies that

$$nE\{(\widehat{\kappa}_1 - \kappa_1)(\widehat{\kappa}_2 - \kappa_2)\} = E[\{\mathcal{F}_1(\bullet) - \kappa_1\}\{\mathcal{F}_2(\bullet) - \kappa_2\}].$$

### A.3. Sketch of Lemma 2.1

We have that

$$\widehat{\kappa}_1 = n^{-1}\sum_{i=1}^{n}\{\delta_i\mathcal{G}(\widetilde{Y}_i) + (1 - \delta_i)\widehat{\mathcal{F}}_i(\bullet)\} = A_1 + A_2,$$

say. By a Taylor series expansion,

$$\begin{aligned}
A_2 &= n^{-1}\sum_{i=1}^{n}(1 - \delta_i)\widehat{\mathcal{F}}_i(\bullet) \\
&= n^{-1}\sum_{i=1}^{n}(1 - \delta_i)[\mathcal{F}_i(\bullet) \\
&\quad + \{\mathcal{F}_{i\mathcal{B}}(\bullet) + \sum_{j=1}^{m}\mathcal{F}_{ij\theta}(\bullet)\theta_{\mathcal{B}}(Z_{ij}, \mathcal{B}_0)\}^{\text{T}}(\widehat{\mathcal{B}} - \mathcal{B}_0) \\
&\quad + \sum_{j=1}^{m}\mathcal{F}_{ij\theta}(\bullet)\{\widehat{\theta}(Z_{ij}, \mathcal{B}_0) - \theta_0(Z_{ij})\}] + o_p(1) \\
&= n^{-1}\sum_{i=1}^{n}(1 - \delta_i)\mathcal{F}_i(\bullet) + \mathcal{M}_{2,\text{imp}}^{\text{T}}(\widehat{\mathcal{B}} - \mathcal{B}_0) \\
&\quad + n^{-1}\sum_{i=1}^{n}\sum_{j=1}^{m}(1 - \delta_i)\mathcal{F}_{ij\theta}(\bullet)\{\widehat{\theta}(Z_{ij}, \mathcal{B}_0) - \theta_0(Z_{ij})\} + o_p(1).
\end{aligned}$$

Using (A.2) and $nh^4 \to 0$, we see

$$n^{-1}\sum_{i=1}^{n}\sum_{j=1}^{m}(1-\delta_i)\mathcal{F}_{ij\theta}(\bullet)\{\widehat{\theta}(Z_{ij},\mathcal{B}_0)-\theta_0(Z_{ij})\}$$

$$= -n^{-2}\sum_{i=1}^{n}\sum_{j=1}^{m}(1-\delta_i)\mathcal{F}_{ij\theta}(\bullet)\sum_{k=1}^{n}\sum_{\ell=1}^{m}\delta_k K_h(Z_{k\ell}-Z_{ij})\mathcal{L}_{k\ell\theta}(\bullet)/\Omega(Z_{ij})$$

$$+n^{-2}\sum_{i=1}^{n}\sum_{j=1}^{m}(1-\delta_i)\mathcal{F}_{ij\theta}(\bullet)\sum_{k=1}^{n}\sum_{\ell=1}^{m}\delta_k\mathcal{L}_{k\ell\theta}(\bullet)\mathcal{G}(Z_{ij},Z_{k\ell})/\Omega(Z_{ij})+o_p(n^{-1/2})$$

$$= n^{-1}\sum_{k=1}^{n}\sum_{\ell=1}^{m}\delta_k\mathcal{L}_{k\ell\theta}(\bullet)C_{1,\mathrm{imp}}(Z_{k\ell})$$

$$+n^{-1}\sum_{k=1}^{n}\sum_{\ell=1}^{m}\delta_k\mathcal{L}_{k\ell\theta}(\bullet)C_{2,\mathrm{imp}}(Z_{k\ell})+o_p(n^{-1/2})$$

$$= n^{-1}\sum_{i=1}^{n}\delta_i D_{i,\mathrm{imp}}(\bullet)+o_p(n^{-1/2}).$$

The result now follows directly from (A.3).

### *A.4. Sketch of Lemma 2.2*

We have that

$$\widehat{\kappa}_2 = n^{-1}\sum_{i=1}^{n}\left[\frac{\delta_i}{\pi(\widetilde{X}_i,\widetilde{Z}_i,\widehat{\zeta})}\mathcal{G}(\widetilde{Y}_i)+\left\{1-\frac{\delta_i}{\pi(\widetilde{X}_i,\widetilde{Z}_i,\widehat{\zeta})}\right\}\widehat{\mathcal{F}}_i(\bullet)\right] = A_1+A_2,$$

say. By a Taylor series expansion,

$$A_1 = n^{-1}\sum_{i=1}^{n}\frac{\delta_i\mathcal{G}(\widetilde{Y}_i)}{\pi(\widetilde{X}_i,\widetilde{Z}_i,\zeta)}$$
$$-E\left\{\frac{\mathcal{G}(Y)}{\pi(\widetilde{X},\widetilde{Z},\zeta)}\pi_\zeta(\widetilde{X},\widetilde{Z},\zeta)\right\}^{\mathrm{T}}n^{-1}\sum_{i=1}^{n}\psi_{i\zeta}+o_p(n^{-1/2}).$$

In addition,

$$A_2 = B_1+B_2+o_p(n^{-1/2});$$
$$B_1 = n^{-1}\sum_{i=1}^{n}\left\{1-\frac{\delta_i}{\pi(\widetilde{X}_i,\widetilde{Z}_i,\zeta)}\right\}\widehat{\mathcal{F}}_i(\bullet);$$
$$B_2 = n^{-1}\sum_{i=1}^{n}\frac{\delta_i\widehat{\mathcal{F}}_i(\bullet)}{\{\pi(\widetilde{X}_i,\widetilde{Z}_i,\zeta)\}^2}\pi_\zeta(\widetilde{X}_i,\widetilde{Z}_i,\zeta)^{\mathrm{T}}(\widehat{\zeta}-\zeta)+o_p(n^{-1/2}).$$

It follows easily that

$$\begin{aligned} B_1 &= n^{-1} \sum_{i=1}^{n} \left\{1 - \frac{\delta_i}{\pi(\widetilde{X}_i, \widetilde{Z}_i, \zeta)}\right\} \mathcal{F}_i(\bullet) \\ &+ n^{-1} \sum_{i=1}^{n} \sum_{j=1}^{m} \left\{1 - \frac{\delta_i}{\pi(\widetilde{X}_i, \widetilde{Z}_i, \zeta)}\right\} \mathcal{F}_{ij\theta}(\bullet)\{\widehat{\theta}(Z_{ij}, \mathcal{B}_0) - \theta_0(Z_{ij})\} \\ &+ n^{-1} \sum_{i=1}^{n} \left\{1 - \frac{\delta_i}{\pi(\widetilde{X}_i, \widetilde{Z}_i, \zeta)}\right\} \mathcal{F}_{i\mathcal{B}}^{\mathrm{T}}(\bullet)(\widehat{\mathcal{B}} - \mathcal{B}_0) \\ &+ n^{-1} \sum_{i=1}^{n} \left\{1 - \frac{\delta_i}{\pi(\widetilde{X}_i, \widetilde{Z}_i, \zeta)}\right\} \sum_{j=1}^{m} \mathcal{F}_{ij\theta}(\bullet) \theta_{\mathcal{B}}^{\mathrm{T}}(Z_i, \mathcal{B}_0)(\widehat{\mathcal{B}} - \mathcal{B}_0) \\ &+ o_p(n^{-1/2}). \end{aligned}$$

Using $0 = E\{1 - \frac{\delta}{\pi(\widetilde{X}, \widetilde{Z}, \zeta)}|\widetilde{X}, \widetilde{Z}\}$ and $\widehat{\mathcal{B}} - \mathcal{B}_0 = O_p(n^{-1/2})$, we see that the third and fourth terms are $o_p(n^{-1/2})$. Also, using $nh^4 \to 0$ and (A.2) we see

$$\begin{aligned} &n^{-1} \sum_{i=1}^{n} \sum_{j=1}^{m} \left\{1 - \frac{\delta_i}{\pi(\widetilde{X}_i, \widetilde{Z}_i, \zeta)}\right\} \mathcal{F}_{ij\theta}(\bullet)\{\widehat{\theta}(Z_{ij}, \mathcal{B}_0) - \theta_0(Z_{ij})\} \\ &= -n^{-1} \sum_{i=1}^{n} \sum_{j=1}^{m} \left\{1 - \frac{\delta_i}{\pi(\widetilde{X}_i, \widetilde{Z}_i, \zeta)}\right\} \mathcal{F}_{ij\theta}(\bullet) \\ &\quad \times \left\{n^{-1} \sum_{k=1}^{n} \sum_{\ell=1}^{m} \delta_k K_h(Z_{k\ell} - Z_{ij}) \mathcal{L}_{k\ell\theta}(\bullet)/\Omega(Z_{ij})\right\} \\ &\quad - n^{-1} \sum_{i=1}^{n} \sum_{j=1}^{m} \left\{1 - \frac{\delta_i}{\pi(\widetilde{X}_i, \widetilde{Z}_i, \zeta)}\right\} \mathcal{F}_{ij\theta}(\bullet) \\ &\quad \times \left\{n^{-1} \sum_{k=1}^{n} \sum_{\ell=1}^{m} \delta_k \mathcal{L}_{k\ell\theta}(\bullet) \mathcal{G}(Z_{ij}, Z_{k\ell})/\Omega(Z_{ij})\right\} + o_p(n^{-1/2}) \\ &= o_p(n^{-1/2}). \end{aligned}$$

The last step follows because $E[\{1 - \frac{\delta}{\pi(\widetilde{X}, \widetilde{Z}, \zeta)}\} \mathcal{F}_{j\theta}(\bullet)|\widetilde{Z}] = 0$. Hence we see that

$$B_1 = n^{-1} \sum_{i=1}^{n} \left\{1 - \frac{\delta_i}{\pi(\widetilde{X}_i, \widetilde{Z}_i, \zeta)}\right\} \mathcal{F}_i(\bullet) + o_p(n^{-1/2}).$$

Similarly, using $\widehat{\zeta} - \zeta = O_p(n^{-1/2})$, $\widehat{\mathcal{B}} - \mathcal{B}_0 = O_p(n^{-1/2})$ and $\widehat{\theta}(z, \mathcal{B}_0) - \theta_0(z) = O_p(n^{-1/2})$, it also follows that

$$B_2 = E\left\{\frac{1}{\pi(\widetilde{X}, \widetilde{Z}, \zeta)} \mathcal{F}(\bullet) \pi_\zeta(\widetilde{X}, \widetilde{Z}, \zeta)\right\}^{\mathrm{T}} n^{-1} \sum_{i=1}^{n} \psi_{i\zeta}(\bullet) + o_p(n^{-1/2}).$$

The result now follows by collecting the terms.

**Acknowledgments.** The authors thank Xihong Lin for providing the data used in this paper.

# References

[1] Chen, X., Linton, O. and Van Keilegom, I. (2003). Estimation of semiparametric models when the criterion function is not smooth. *Econometrica* **71** 1591–1608. MR2000259

[2] Lin, X. and Carroll, R. J. (2006). Semiparametric estimation in general repeated measures problems. *J. R. Stat. Soc. Ser. B Stat. Methodol.* **68** 69–88. MR2212575

[3] Maity, A., Ma, Y. and Carroll, R. J. (2007). Efficient estimation of population-level summaries in general semiparametric regression models. *J. Amer. Statist. Assoc.* **102** 123–139. MR2293305

[4] Sepanski, J. H., Knickerbocker, R. and Carroll, R. J. (1994). A semiparametric correction for attenuation. *J. Amer. Statist. Assoc.* **89** 1366–1373. MR1310227

[5] Zeger, S. and Diggle, P. (1994). Semiparametric models for longitudinal data with application to cd4 cell numbers in hiv seroconverters. *Biometrics* **50** 689–699.

[6] Zhang, D., Lin, X., Raz, J. and Sowers, M. (1998). Semiparametric stochastic mixed models for longitudinal data. *J. Amer. Statist. Assoc.* **93** 710–719. MR1631369